\documentclass{amsart}%
\usepackage{amsfonts}
\usepackage{amsmath}
\usepackage{amssymb}
\usepackage{graphicx}%
\providecommand{\U}[1]{\protect\rule{.1in}{.1in}}
\newtheorem{theorem}{Theorem}
\theoremstyle{plain}

\newtheorem{definition}{Definition}

\newtheorem{lemma}{Lemma}

\newtheorem{remark}{Remark}

\numberwithin{equation}{section}
\begin{document}
\title[Finiteness of Subfamilies of Calabi-Yau n-Folds over Curves]{Finiteness of Subfamilies of Calabi-Yau n-Folds over Curves with Maximal
Length of Yukawa-Coupling}
\author{Kefeng Liu}
\address{UCLA, Department of Mathematics\\
Los Angeles, CA 90095\\
and Center of Mathematical Sciences\\
Zheijiang University}
\author{Andrey Todorov}
\address{UC, Department of Mathematics\\
Santa Cruz, CA 95064\\
Institute of Mathematics\\
Bulgarian Academy of Sciences}
\author{Shing-Tung Yau}
\address{Harvard University\\
Department of Mathematics\\
Cambridge, Mass. 02138}
\author{Kang Zuo}
\address{University Mainz\\
Department of Mathematics\\
Mainz 55099    } \dedicatory{Dedicated to the memory of E. Viehweg}

\footnotetext[1]{The last named author  was supported by SFB/Transregio 45
Periods, Moduli Spaces and Arithmetic of Algebraic Varieties of the
DFG (Deutsche Forschungsgemeinschaft).}

\begin{abstract}

We give a simplified and more algebraic proof of the finiteness of
the families of Calabi-Yau n-folds with non-vanishing of
Yukawa-coupling over a fixed base curve and with fixed degeneration
locus. We also give a generalization of this result. Our method is
variation of Hodge structure and poly-stability of Higgs bundles.

\end{abstract}
\maketitle

\section{\textbf{Introduction}}

In this short note, we give a simplified  proof of the main theorem
in \cite{LTYZ03}, namely the set of subfamilies of Calabi-Yau
n-folds with non-vanishing of Yukawa-coupling over the fixed base
curve and with fixed degeneration locus is finite. The method is
based on variation of Hodge structure and poly-stability of Higgs
bundles. We hope this approach will be more understandable to
algebraic geometers, particularly to people working on Hodge theory.
Note that if we consider the subset of sub-families passing through
the maximal nilpotent degeneration locus, the finiteness of such
families was proved by Zhang Yi in his PhD thesis \cite{Zh04}.
However, in general a rigid family does not need passing
through the maximal nilpotent degeneration locus. \\[0.1cm]Very recently J. C.
Rohde \cite{rohde09} and later A. Garbagnati and B. Van Geemen
\cite{Gar-Gee09} have found some universal families of Calabi-Yau
3-folds without maximal nilpotent degeneration locus. Moreover the
Yukawa-coupling vanishes.  To cover more general cases, in this note
we also formulate a more general criterion for the rigidity of
sub-families of those universal families with vanishing
Yukawa-coupling. We use the notion of the length of Yukawa coupling
attached to a subfamily introduced in \cite{VZ03} and show that the
set of subfamilies of Calabi-Yau n-folds over the fixed
base curve and with fixed degeneration locus, whose Yukawa-coupling have
the same length as the length of the universal family, is finite.\\[.2cm]

{\bf Acknowledgements:}  This work was done while the last named author was visiting
Center of Mathematical Sciences at Zhejing University and East China Normal University. He would like
to thank both institutions' financial support and the hospitality.\\[.2cm]

\section{\textbf{Universal Family over Moduli Space with Level $N-$Structure}}

Let ${\mathcal{M}}_{h}({\mathbb{C}})$ denote the set of isomorphic classes of
minimal polarized manifolds $F$ with fixed Hilbert polynomial $h$, and let
${\mathcal{M}}_{h}$ be the corresponding moduli functor, i.e.
\[
{\mathcal{M}}_{h}(U)=\left\{
\begin{array}
[c]{c}%
(f:V\rightarrow U,{\mathcal{L}});f\mbox{ smooth and}\\
(f^{-1}(u),{\mathcal{L}}|_{f^{-1}}(u))\in{\mathcal{M}}_{h}({\mathbb{C}%
}),\mbox{ for all }u\in U
\end{array}
\right\}
\]
By Viehweg \cite{v95} there exists a quasi-projective coarse moduli scheme
$M_{h}$ for ${\mathcal{M}}_{h}$. Fixing a projective manifold $\bar{U}$ and
the complement $U$ of a normal crossing divisor, we want to consider
\[
\mathrm{\mathbf{H_{\bar{U},U}}}=\left\{
\begin{array}
[c]{c}%
\varphi:(\bar{U},U)\rightarrow(\bar{M}_{h},M_{h})\quad\mathrm{induced}\\
\mathrm{by\,non-isotrivial\, families}\,f:X\rightarrow U
\end{array}
\right\}  .
\]
Since $M_{h}$ is just a coarse moduli scheme, it is not clear whether
$\mathrm{\mathbf{H}}$ has a scheme structure. However, since we are
considering here the moduli space of polarized Calabi-Yau n-folds, all
$F\in{\mathcal{M}}({\mathbb{C}})$ admit locally injective Torelli map for the
Hodge structure on $H^{n}(F,\mathbb{Z}).$ Hence by Theorem 10.3, f) in
\cite{Popp77}(see also \cite{Sz98}) there exists a fine moduli scheme
$M_{h}^{N}$ with a sufficiently large level structure $N$ on $H^{n}%
(F,\mathbb{Z}),$ which is finite over $M_{h}.$ By abuse of notations, we will
replace ${\mathcal{M}}_{h}$ by the moduli functor of polarized manifolds with
a level $N$ structure, and choose some smooth compactification $\bar
{M}_{h}$, such that $\bar{M_{h}}\setminus M_{h}=:S$ is a normal crossing
divisor.. Then $\mathrm{\mathbf{H_{\bar{U},U}}}$ parameterizes all morphisms
from $\varphi:(\bar{U},U)\rightarrow(\bar{M}_{h},M_{h})$, hence it is a
scheme. Moreover there exists a universal family $f:X\rightarrow
\mathrm{\mathbf{H_{\bar{U},U}}}\times U.$

\section{\textbf{Boundedness of Families of Polarized Manifolds}}

\begin{theorem}
\label{1} Assume $\dim\bar{U}=1,$ then $\mathrm{\mathbf{H_{\bar{U},U}}}$ is of
finite type.
\end{theorem}

Theorem \ref{1} is known since long time for many cases and due to
many people: the case of curves is due to Arakelov and Parshin, and
the case of abelian varieties is due to Faltings \cite{fa83}. Yau
\cite{ya78} proved a generalization of Schwarz lemma which implies
the Arakelov inequalities and also that Arakelov equality holds for
a family of curves of genus bigger or equal to two if and only if
the VHS split. This result was obtained earlier by Viehweg and K.
Zuo using the technique of Higgs bundles. The Schwarz lemma of Yau
was used by Jost-Yau \cite{jy93} to give another proof of Falting's
theorem in \cite{fa83}. C. Peters \cite{Pe90} has proven the
boundedness for general period maps. Recently this kind of
boundedness result has been generalized to polarized manifolds
without assuming the local injectivity of Torelli map. The case of
surfaces of general type is due to Bedulev-Viehweg \cite{bv00}, the
case of surfaces of Kodaira dimension $\leq1$ is due to
Oguiso-Viehweg \cite{ov01} and the case the semi positive canonical
line bundle is due do Viehweg-Zuo \cite{vz02}. Very recently, Kovace
and Lieblich \cite{kl06} have shown the boundedness for canonically
polarized manifolds over higher dimensional bases.

As the case of Calabi-Yau manifolds is a sub-case of the
semi-positive canonical line bundle, for reader's convenience we
reproduce the proof for
Theorem 1, which can be found in (\cite{vz02}, Theorem 6.2).\\[0.1cm]%
\textbf{Proof of Theorem 1.}:Viehweg \cite{v95} has constructed a
series of nef invertible sheaves $\lambda_{\nu}$ on $M_{h}$ such
that $\lambda_{\nu}$ are ample $\nu>>1$ and are natural in the sense
$\varphi^{\ast}\lambda_{\nu }=\det f_{\ast}\omega_{X/U}^{\nu}$ for
$\varphi:U\rightarrow M_{h}$ induced by a family $f:X\rightarrow U.$
He has shown furthermore \cite{v06} that $\lambda_{\nu}$ extends to
nef invertible sheaves $\bar{\lambda}_{\nu}$ on $\bar{M}_{h}$ which
are again natural in the sense that if $\varphi:\bar
{U}\rightarrow\bar{M}_{h}$ is induced by a semi stable family
$\bar {f}:\bar{X}\rightarrow\bar{U}$ then
\[
\varphi^{\ast}\bar{\lambda}_{\nu}=\det\bar{f}_{\ast}\omega_{\bar{X}/\bar{U}%
}^{\nu}.
\]

For a point $(\varphi:(\bar{U},U)\rightarrow(\bar{M}_{h},M_{h}%
))\in\mathrm{\mathbf{H_{\bar{U},U}}},$ let $f:X\rightarrow U$ denote
the pull-back of the universal family via $\varphi.$ Assume
$\dim\bar{U}=1$ and let $S:=\bar{U}\setminus U$, $s:=\#S$, $F$ the
generic fibre of $f$ and $n:=\dim F$ then one finds the Arakelov
inequality for $f$ \cite{vz01} and \cite{vz02}
\[
\deg(\varphi^{\ast}\bar{\lambda}_{\nu})\leq\nu\cdot\mathrm{rank}f_{\ast}%
\omega_{X/U}^{\nu}\cdot(n\cdot\deg\Omega_{\bar{U}}^{1}(\log S)+s)\cdot
e(\omega_{F}^{\nu}),
\]
where $e(\omega_{F}^{\nu}))$ is the holomorphic Euler characteristic of
$\omega_{F}^{\nu}$. Fixing a $\nu>>1,$ then $\bar{\lambda}_{\nu}$ is nef on
$\bar{M}_{h}$, ample on $M_{h}$ and $\deg(\varphi^{\ast}\bar{\lambda
}_{\nu})$ is particularly bounded above by a number independent of $\varphi$.
Hence $\mathrm{\mathbf{H_{\bar{U},U}}}$ is of finite type. Theorem 1 is proved.

\section{\textbf{Higgs Bundles and Yukawa-Couplings Attached to Families,
Rigidity and Finiteness}}

Let $\mathbb{V}$ denote the local system
$R^{n}f_{univ\ast}\mathbb{Z}_{X}$ attached to the universal family
$f_{univ}:X\rightarrow M_{h}.$ Then the local monodromies around each
component of $S$ are quasi-unipotent. By taking another finite cover
of $\bar{M_{h}}$, which is \'{e}tale over $M_{h}$ we may assume
the local monodromies around all components of $S$ are unipotent.
$\mathbb{V}$ carries a polarized $\mathbb{Z}-$variation of Hodge
structure of weight-$n$ with the Hodge filtration
\[
F^{n}\subset\cdots F^{0}=\mathbb{V}\otimes\mathcal{O}_{U}=:V,
\]
satisfying the Griffiths's transversality. The grading $F^{p}/F^{p+1}%
=:E^{p,n-p}$ with the projection Gauss-Manin connection $\theta^{p,q}%
:E^{p,q}\rightarrow E^{p-1,q+1}\otimes\Omega_{M_{h}}^{1}$ defines the
correspondence between $\mathbb{V}$ and its Higgs bundle
\[
(E:=\bigoplus_{p+q=n}E^{p,q},\theta:=\bigoplus_{p+q=n}\theta^{p,q}%
),\quad\theta\wedge\theta=0.
\]
The following properties are known
\[
E^{p,q}\simeq R^{q}f_{univ\ast}\Omega_{X/M_{h}}^{p},
\]
and the Higgs map
\[
\theta^{p,q}:R^{q}f_{univ\ast}\Omega_{X/M_{h}}^{p}\rightarrow R^{q+1}f_{univ\ast
}\Omega_{X/M_{h}}^{p-1}\otimes\Omega_{M_{h}}^{1}%
\]
coincides with the Kodaria-Spencer maps of $f_{univ}$ on $R^{q}f_{univ\ast}%
\Omega_{X/M_{h}}^{p}$. By Deligne $(E,\theta)$ admits a canonical extension,
such that the extended $\theta$ takes the value in $\Omega_{\bar{M_{h}}%
}^{1}(\log S).$ The extension is denoted again by $(E,\theta).$

For a point $(\varphi:(\bar U, U) \to(\bar M_{h}, M_{h}))\in
\mathrm{\mathbf{H_{\bar U, U}}}$, let $f: X\to U$ denote the
pull-back of the universal family via $\varphi.$ Then the pull-back
$\varphi ^{*}(\mathbb{V})$ is the VHS $\mathbb{V}_{f}$ attached to
$f$. Let $(E_{f},\theta_{f})$ be the Higgs bundle corresponding to
$\mathbb{V}_{f},$ then $E_{f}=\varphi^{*}E,$ and $\theta_{f}$ is
equal to $d\varphi(\theta)$ defined as
\[
d\varphi(\theta): E_{f}\overset{\varphi^{*}\theta}{\rightarrow}E_{f}%
\otimes\varphi^{*}\Omega^{1}_{\bar M_{h}}(\log S)\overset{d\varphi
}{\rightarrow}E_{f}\otimes\Omega^{1}_{\bar U}(\log S)
\]
\\[.1cm]

\begin{definition}
The $n$-times iterated Kodaira-Spencer map
\[
\theta_{f}^{n}:E_{f}^{n,0}\overset{\theta_{f}^{n,0}}{\rightarrow}E_{f}%
^{n-1,1}\otimes\Omega_{U}^{1}\overset{\theta_{f}^{n-1,1}}{\rightarrow}%
\cdots\overset{\theta_{f}^{1,n-1}}{\rightarrow}E_{f}^{0,n}\otimes S^{n}%
\Omega_{\bar{U}}^{1}(\log S)
\]
is called the Yukawa-coupling $\theta_{f}^{n}$ attached to the family
$f:X\rightarrow U.$ We will say that the Yukawa-coupling attached the $f$ does
not vanish, if $\theta_{f}^{n}\not =0.$
\end{definition}

In general, following  \cite{VZ03} we consider the $l$-times iterated
Kodaira-Spencer map attached to the family $f:X\to U$ for $1\leq l\leq n.$%

\[
\theta^{l}_{f}: E^{n,0}_{f} \overset{\theta^{n,0}_{f}}{\rightarrow}%
E^{n-1,1}_{f}\otimes\Omega^{1}_{U}\overset{\theta^{n-1,1}_{f}}{\rightarrow
}\cdots\overset{\theta^{n-l+1,l-1}_{f}}{\rightarrow} E^{n-l,l}_{f}\otimes
S^{l}\Omega^{1}_{\bar U}(\log S).
\]
\\[.2cm]

\begin{definition}
The length of the Yukawa-coupling attached to a family $f:X\rightarrow U$ is
the maximal $l$ such that $\theta_{f}^{l}\not =0.$ This number is an invariant
of the family $f:X\rightarrow U$ and will be denoted by $l(\theta_f)$.
\end{definition}

\begin{remark}
The length of the Yuwawa-coupling attached to a sub family can be in general
smaller than the Yuwawa-coupling attached to the universal family.
\\[0.1cm]The length of the universal family of Calabi-Yau quintic hyper
surfaces in $P^{4}$ is equal to 3. By  \cite{VZ03} there do exist
sub-families of Calabi-Yau quintic hypersurfaces with the length of
Yukawa-coupling equal to 1, or 2.\\
Very recently Rohde \cite{rohde09} has found a universal family of Calabi-Yau
3-folds over a complex ball quotient, which does not have maximal nilpotent
degeneration points. (see also \cite{Gar-Gee09}). Moreover the length of the
Yukawa-coupling of this universal family is equal to 1.
\end{remark}

\begin{definition}
We will say that the family $f:X\rightarrow U$ has the Yukawa-coupling with
the maximal length if $l(\theta_{f})=l(\theta_{f_{univ}}).$
\end{definition}

\begin{theorem}
\label{2.a} Assume $\theta_{f}^{n}$ attached to a family $f:X\rightarrow U$
does not vanish, then $f$ is rigid, i.e. the component of
$\mathrm{\mathbf{H_{\bar{U},U}}}$ containing $\varphi:(\bar{U},U)\rightarrow
(\bar{M}_{h},M_{h})$ induced by $f:X\rightarrow U$ is zero-dimensional.
\end{theorem}

Regarding the existence of universal families of Calabi-Yau $n$-folds whose
Yukawa-coupling does vanish, i.e. the length of the Yukawa-coupling is smaller
than $n$ we have the following corresponding criterion for the rigidity.

\begin{theorem}
\label{2.b}Assume that the length of the Yukawa-coupling attached to a family
$f:X\rightarrow U$ is equal to the length of the Yukawa-coupling attached the
universal family, then $f$ is rigid.
\end{theorem}

Theorem \ref{2.a} and \ref{2.b} combined together with Theorem \ref{1} on the
boundedness of $\mathrm{\mathbf{H_{\bar{U},U}}}$ imply

\begin{theorem}
\label{3.a}Assume $\dim U=1,$ then the subset
\[
\mathrm{\mathbf{H_{\bar{U},U}}}^{0}:=\{(\varphi:(\bar{U},U)\rightarrow
(\bar{M}_{h},M_{h}))\in\mathrm{\mathbf{H_{\bar{U},U}}}\,|\,\theta_{f}%
^{n}\not =0\}\subset\mathrm{\mathbf{H_{\bar{U},U}}}%
\]
is finite.
\end{theorem}

Theorem \ref{3.a} can be generalized as follows:

\begin{theorem}
\label{3.b}Assume $\dim U=1,$ then the subset%
\[
\mathrm{\mathbf{H_{\bar{U},U}}}^{0}:=\{(\varphi:(\bar{U},U)\rightarrow
(\bar{M}_{h},M_{h}))\in\mathrm{\mathbf{H_{\bar{U},U}}}\,|\,l(\theta
_{f})=l(\theta_{f_{univ}})\}\subset\mathrm{\mathbf{H_{\bar{U},U}}}%
\]
is finite.
\end{theorem}
 As a corollary from  Theorem 5 we have

\begin{theorem}
If $ l(\theta_{f_{univ}})=1$  (for example, J. C. Rohde's
universal family) then $\mathrm{\mathbf{H_{\bar{U},U}}}$ is
finite.
\end{theorem}

\section{\textbf{Decomposition of $(E_{f},\theta_{f})$ Induced by
Infinitesimal Deformations and Proofs of Theorems \ref{2.a}, \ref{2.b},
\ref{3.a} and \ref{3.b}}}

The universal logarithmic Higgs bundle
\[
\theta: E\to E\otimes\Omega^{1}_{\bar{M_{h}}}(\log S)
\]
attached to the universal family $f_{univ}:X\to M_{h}$ tautologically defines the
universal Higgs map
\[
\theta: T_{\bar{M_{h}}}(-\log S)\to End(E),
\]
which coincides with the differential of the period map of the universal
family. By the Griffiths's transversality $\theta$ is of the Hodge type
(-1,1), i.e. for any local vector field $v$ of $T_{\bar{M_{h}}}(-\log S)$
then $\theta_{v}(E^{p,q})\subset E^{p-1,q+1}.$

A non-trivial infinitesimal deformation of $\varphi:(\bar U, U) \to(\bar
M_{h}, M_{h})$ is a non-zero section
\[
\mathcal{O}_{\bar U}\overset{t}{\rightarrow}\varphi^{*} T_{\bar{M_{h}}%
}(-\log S),
\]
which does not factor through $d\varphi: T_{\bar U}(-\log S)\to\varphi^{*}
T_{\bar{M_{h}}}.$

The composition of $t$ with the pull-back of the universal Higgs map
\[
\varphi^{\ast}\theta_{t}:\mathcal{O}_{\bar{U}}\overset{t}{\rightarrow}%
\varphi^{\ast}T_{\bar{M_{h}}}(-\log S)\overset{\varphi^{\ast}%
\theta}{\rightarrow}End(E_{f})
\]
corresponds to the induced infinitesimal deformation of the period map of the
family $f:X\rightarrow U$.

\begin{lemma}
\label{5}\textbf{i)} $\varphi^{\ast}\theta_{t}$ is an endomorphism of the
Higgs bundle $(E_{f},\theta_{f}).$\\[0.1cm]\textbf{ii)} Let $K:=\mathrm{ker}%
(\varphi^{\ast}\theta_{t})\subset E_{f}.$ Then $K$ is a Higgs subbundle and
there exists a decomposition of the Higgs bundle $(E_{f},\theta_{f})$
\[
(E_{f},\theta_{f})=(K,\theta_{f}|_{K})\oplus(K^{\bot},\theta_{f}|_{K^{\bot}}),
\]
where $K^{\bot}$ is orthogonal to $K$ w.r.t the Hodge metric.
\end{lemma}

\textbf{Proof of Lemma} \ref{5}: We give two proofs for i). First proof of
\textbf{i)}: By Theorem 3.2, (i) in \cite{Pe90} $\varphi^{\ast}\theta_{t}$ is
induced by a flat endomorphism
\[
\alpha:\mathbb{V}_{f}\otimes\mathbb{C}\rightarrow\mathbb{V}_{f}\otimes
\mathbb{C}%
\]
of the Hodge type (-1,1), i.e. $\alpha\in End(V_{f})$, commutes with the
Gauss-Manin connection, and $\alpha(F_{f}^{p})\subset F_{f}^{p-1}.$ The
induced endomorphism by $\alpha$ on the gradings $E=\bigoplus_{p+q=n}%
E_{f}^{p,q}$ coincides with $\varphi^{\ast}\theta_{t}.$ Note that the Higgs
map $\theta_{f}$ is defined by the projection of the Gauss-Manin connection of
the gradings, the commutativity of $\alpha$ with the Gauss-Manin connection
implies the commutativity of $\varphi^{\ast}\theta_{t}$ with $\theta_{f},$
i.e., $\varphi^{\ast}\theta_{t}$ is an endomorphism of Higgs bundles
\[
\varphi^{\ast}\theta_{t}:(E_{f},\theta_{f})\rightarrow(E_{f},\theta_{f}).
\]
\\[0.1cm]The second proof of \textbf{i)} goes back to \cite{z00}, and is much
simpler. Consider the universal Higgs map on $\bar{M_{h}}$
\[
\theta:E\rightarrow E\otimes\Omega_{\bar{M}_{h}}^{1}(\log S),\quad
\quad\theta\wedge\theta=0\,\,\mathrm{in}\,\,End(E)\otimes\Omega_{\bar
{M}_{h}}^{2}(\log S).
\]
If $e_{1},\cdots,e_{m}$ is a local basis of $\Omega_{\bar{M}_{h}}%
^{1}(\log S)|_{W}$ over some open set $W$, then
\[
\theta=\sum\theta_{i}e_{i},
\]
where $\theta_{i}$ are endomorphisms of $E|_{W}$ The condition $\theta
\wedge\theta=0$ means that $\theta_{i}$ commutes with one another. Given a
local vector field
\[
v=\sum c_{i}e_{i}^{\ast}\in T_{\bar{M}_{h}}(-\log S)|_{W},
\]
the image of $v$ under the map
\[
\theta:T_{\bar{M}_{h}}(-\log S)|_{W}\rightarrow End(E|_{W})
\]
defines an endomorphism
\[
\theta_{v}=\sum c_{i}\theta_{i}:E_{W}\rightarrow E_{W},
\]
then $\theta_{v}$ commutes with
\[
\theta:E_{W}\rightarrow E_{W}\otimes\Omega_{\bar{M}_{h}}^{1}(\log
S)|_{W}.
\]
Pulling back this commutativity via $\varphi,$ we see that $\varphi^{\ast
}\theta_{t}$ commutes with $\varphi^{\ast}\theta$, hence commutes with
$d\varphi(\theta)=\theta_{f}.$ $\square$

As for \textbf{ii)} By \cite{si92}, Section 4, it is known that the
Hodge metric on $(E_{f},\theta_{f})|_{U}$ is Hermitian-Yang-Mills
and by \cite{si90}, Theorem 5, the canonical extended Higgs bundle,
$(E_{f},\theta_{f})$ is Higgs poly-stable of slope zero. By i)
$\varphi^{\ast}\theta_{t}$ is an endomorphism of
$(E_{f},\theta_{f})$, $\mathrm{ker}(\varphi^{\ast}\theta_{t})(=K)$
and $\mathrm{im}(\varphi^{\ast }\theta_{t})$ are Higgs subsheaves.

Applying the Higgs-poly-stability to $K$ and $\mathrm{im}(\varphi^{\ast}%
\theta_{t})$ one finds
\[
\deg K\leq0,\quad\deg\mathrm{im}(\varphi^{\ast}\theta_{t})\leq0.
\]
Since
\[
\deg K+\deg\mathrm{im}(\varphi^{\ast}\theta_{t})=\deg(E_{f})=0,
\]
$\deg K=0.$ Again by applying the Higgs poly-stability on $K$, $K$
defines a splitting
\[
(E_{f},\theta_{f})=(K,\theta_{f}|_{K})\oplus(K^{\bot},\theta_{f}|_{K^{\bot}%
}),
\]
such that $K^{\bot}$ is orthogonal to $K$ w.r.t. the Hodge metric. ii) is
complete. $\square$

\begin{lemma}
\label{6} \textbf{i)} $E_{f}^{0,n}\subset K.$ \quad\textbf{ii)} $E_{f}%
^{n,0}\subset K^{\bot}.$
\end{lemma}

\textbf{Proof of Lemma }\ref{6}: Since $\varphi^{\ast}\theta_{t}$ is of the
Hodge type (-1,1),
\[
\varphi^{\ast}\theta_{t}(E_{f}^{0,n})\subset E_{f}^{-1,n+1}=0.
\]

i) is proven. $\square$

As for ii). We note first that the Kodaira-Spencer map
\[
\theta: T_{M_{h}}\otimes E^{n,0}\to E^{n-1,1}%
\]
for the universal family of Calabi-Yau manifolds is an isomorphism.
Hence, particularly, the sheaf morphism $\varphi^{*}\theta_{t}:
E^{n,0}_{f}\to
E^{n-1,1}_{f}$ is injective.\\[.1cm]

Write any local section
\[
s\in E_{f}=\bigoplus_{p+q=n}E_{f}^{p,q}%
\]
in the form
\[
s=\sum_{p+q=n}s^{p,q},
\]
then $\varphi^{\ast}\theta_{t}$ maps $s^{p,q}$ into $E_{f}^{p-1,q+1}.$ Hence
if $\varphi^{\ast}\theta_{t}(s)=0,$ then all $\varphi^{\ast}\theta_{t}%
(s^{p,q})=0.$ Thus, the injectivity of $\varphi^{\ast}\theta_{t}$ on
$E_{f}^{n,0}$ implies that $s^{n,0}=0$ if $s\in K,$ i.e.
\[
K\subset E_{f}^{n-1,1}\oplus\cdots\oplus E_{f}^{0,n}.
\]

Since the Hodge decomposition
\[
E_{f}=\bigoplus_{p+q=n}E_{f}^{p,q}%
\]
is orthogonal w.r.t. the Hodge metric, $E_{f}^{n,0}$ is orthogonal
to $K$. ii) is complete.

\textbf{Proof of Theorem} \ref{2.a}: If the component of
$\mathrm{\mathbf{H_{\bar{U},U}}}$ contains a point
\[
\varphi:(\bar{U},U)\rightarrow(\bar {M}_{h},M_{h})
\]
induced by $f:X\rightarrow U$ with $\theta_{f}^{n}\not =0$ is positive
dimensional. Then by Lemma 5 a non-trivial infinitesimal deformation of
$\varphi$ induces a decomposition
\[
(E_{f},\theta_{f})=(K,\theta_{f}|_{K})\oplus(K^{\bot},\theta_{f}|_{K^{\bot}%
}).
\]
By lemma 6
\[
E_{f}^{0,n}\subset K,\quad E_{f}^{n,0}\subset K^{\bot}.
\]
Since $K^{\bot}$ is $\theta_{f}$-invariant,
\[
\theta_{f}^{n}(E^{n,0})\subset K^{\bot}\otimes S^{n}(\Omega_{\bar{U}}^{1}(\log
S).
\]
On the other hand, $\theta_{f}^{n}(E_{f}^{n,0})\subset E_{f}^{0,n}\otimes
S^{n}(\Omega_{\bar{U}}^{1}(\log S)$, which is contained in 
$$K\otimes S^{n}(\Omega_{\bar{U}}^{1}(\log S).$$ 
Since $K\cap K^{\bot}=\{0\},\theta
_{f}^{n}(E^{n,0})=0.$ A contradiction to the assumption $\theta_{f}^{n}%
\not =0.$ The proof of Theorem \ref{2.a} is complete. $\square$

The proof of Theorem \ref{2.b} is similar to the proof of Theorem
\ref{2.a}. The only difference is that Lemma \ref{6} has to be
reformulated in the following form:

We denote
\[
K^{n-l,l}:=\mathrm{Ker}(\theta_{f_{univ}}: E^{n-l,l}\to E^{n-l-1, l+1}%
\otimes\Omega^{1}_{\bar M_{h}}(\log S)).
\]
\\[.2cm]

\begin{lemma}
\label{6.b} \textbf{i)} $\varphi^{\ast}K^{n-l,l}\subset K.$ \quad\textbf{ii)}
$E^{n,0}\subset K^{\bot}.$
\end{lemma}

\textbf{Proof of Lemma }\ref{6.b}: The proof of part ii) of Lemma
\ref{6.b} is exactly the same as the proof of part ii) Lemma
\ref{6}. So we  need to prove only part i). The definition of
$K^{n-l,l}$ implies
\[
\varphi^{\ast}\theta_{f_{univ}}:\varphi^{\ast}T_{\bar{M}_{h}}(-\log
S)\otimes\varphi^{\ast}K^{n-l,l}\rightarrow0
\]
and since $t:\bar{U}\rightarrow\varphi^{\ast}T_{\bar{M}_{h}}(-\log S)$, we
have
\[
\varphi^{\ast}\theta_{t}(\varphi^{\ast}K^{n-l,l})=\varphi^{\ast}\theta_{f_{univ}}(
t\otimes \varphi^{\ast}K^{n-l,l})=0.
\]
This shows i). $\square$

\textbf{Proof of Theorem} \ref{2.b}: Since $l(\theta_{f_{univ}})=l,$ we have
$\theta_{f_{univ}}^{l+1}=0$. This implies
\[
\theta_{f_{univ}}^{l}(S^lT_{\bar{M}_{h}}(-\log S)\otimes E^{n,0})\subset
K^{n-l,l}.
\]
Pulling back this inclusion to $\bar{U}$ via $\varphi:\bar{U}\rightarrow
\bar{M}_{h}$ and by part i) of Lemma \ref{6.b} one obtains
\[
\theta_{f}^{l}(S^lT_{\bar{U}}(-\log S)\otimes E_{f}^{n,0})\subset\varphi^{\ast
}\theta_{f_{univ}}^{l}(S^lT_{\bar{M}_{h}}(-\log S)\otimes E^{n,0})\subset
\varphi^{\ast}K^{n-l,l}\subset K.
\]
On the other hand, since $K^{\bot}$ is $\theta_{f}$-invariant and $E^{n,0}_f\subset K^{\bot}$
\[
\theta_{f}^{i}(S^iT_{\bar{U}}(-\log S)\otimes E^{n,0})\subset K^{\bot}
\]
for $0\leq i\leq n.$ For $i=l$ this inclusion together with the inclusion
above implies that
$$\theta_{f}^{l}(S^lT_{\bar{U}}(-\log S)\otimes E^{n,0})=0$$
which contradicts to $\theta
_{f}^{l}\not =0.$ The proof of Theorem \ref{2.b} is complete. $\square$

\begin{remark}
In the proof of the above theorems on the rigidity, what we need is
the existence of a polarized complex VHS over the moduli space such
that the first Hodge bundle is a line bundle and the first Higgs map
is an isomorphism. Such type VHS is called Calabi-Yau like VHS. In
\cite{VZ03} it is shown that the universal family of hypersurfaces
in $P^{n}$ of degree $\geq n+1$ carries a Calabi-Yau like VHS.
Hence, all theorems here holds true for those subfamilies.
\end{remark}


\begin{thebibliography}{999999999}                                                                                        %


\bibitem[BV00]{bv00}\textsc{E. Bedulev and E. Viehweg:} On the Shafarevich
conjecture for surfaces of general type over function fields. Invent. Math.
139, (2000), 603-615.

\bibitem[De87]{de87}\textsc{P. Deligne:} Un theoreme de finitude pour la
monodromie. Discrete Groups in Geometry and Analysis, Birkh\"auser, Progress
in Math. 67, (1987), 1-19.

\bibitem[Fa83]{fa83}\textsc{G. Faltings:} Arakelov's Theorem for Abelian
Varieties. Invent. Math. 73, (1983), 337-348.

\bibitem[Gar-Gee09]{Gar-Gee09}\textsc{A. Garbagnati and B. Van Geemen:} The
Picard-Fuchs equation of a family of Calabi-Yau three folds without maximal
unipotent monodromy. Math.AG/09034404v1

\bibitem[JY93]{jy93}\textsc{J. Jost and S-T. Yau:} Harmonic Mappings and
Algebraic Varieties Over Function Fields. Am. J. of Math. 115, (1993), 1197-1227.

\bibitem[KL06]{kl06}\textsc{S. Kovacs and M. Lieblich:} Boundedness of
families of canonically polyrized manifolds: a higher dimensional analogue of
Shafarevich conjecture. MathAG/0611672.

\bibitem[LTYZ03]{LTYZ03}\textsc{K. Liu, A. Todorov, S.-T. Yau, and K.
Zuo:}Shafarevich\textexclamdown
\={}%
s conjecture for CY manifolds. I. Q. J. Pure Appl. Math., 1(1):28%
\"{}%
C67, 2005.

\bibitem[OV01]{ov01}\textsc{K. Oguiso and E. Viehweg:} On the isotriviality of
families of elliptic surfaces. J. Alg. Geom. 10, (2001), 569-598.

\bibitem[Popp77]{Popp77}H. Popp.: Moduli theory and classification theory of
algebraic varieties. Lect. Notes in Math. \textbf{620} Springer-Verlag,
Berlin-New York, 1977

\bibitem[Pe90]{Pe90}\textsc{C.A. Peters} Rigidity for variation of Hodge
structure and Arakelov-type finiteness theorems. Compositio Math. 75 (1990) 113-126.

\bibitem[rohde09]{rohde09}\textsc{J.C. Rohde:} Maximal automorphisms of
Calabi-Yau manifolds versus maximally unipotent monodromy. MathAG/0902.4592v2

\bibitem[Si90]{si90}\textsc{C. Simpson:} Harmonic bundles on non compact
curves J. of AMS 3, (1990), 713-770.

\bibitem[Si90]{si92}\textsc{C. Simpson:} Higgs bundles and local systems. IHES
75 (1992) 5-95.

\bibitem[Sz98]{Sz98}\textsc{B. Szendr{o}i} Some finiteness results for
Calabi-Yau threefolds. ArXiv: alg-geom/9708011v3.

\bibitem[V95]{v95}\textsc{E. Viehweg:} Quasi-projective Moduli for Polarized
Manifolds. Ergebnisse der Mathematik, 3. Folge 30 (1995), Springer Verlag,
Berlin-Heidelberg-New York.

\bibitem[V06]{v06}\textsc{E. Viehweg:} Compactifications of smooth Families
and of moduli spaces of polarized manifolds. Math.AG/0605093

\bibitem[VZ01]{vz01}\textsc{E. Viehweg and K. Zuo:} On the isotriviality of
families of projective manifolds over curves. J. of Alg. Geom. 10, (2001), 781-799.

\bibitem[VZ02]{vz02}\textsc{E. Viehweg and K. Zuo:} Base spaces of
non-isotrivial families of smooth minimal models. In: Complex Geometry
(Collection of Papers dedicated to Hans Grauert), Springer Verlag,
Berlin-Heidelberg-New York (2002), 279-328.

\bibitem[VZ03]{VZ03}\textsc{E. Viehweg and K. Zuo:} Complex multiplication,
Griffiths-Yukawa couplings, and rigidity for families of hypersurfaces. J.
Algebraic Geom. 14 (2005), no. 3, 481--528.

\bibitem[Yau78]{ya78}\textsc{S-T. Yau:} A general Schwarz lemma for K\"ahler
manifolds. Am. J. of Math. 100 (1978),197-203.

\bibitem[Zh04]{Zh04}\textsc{Y. Zhang:} Rigidity for families of polarized
Calabi-Yau varieties. J. Differential Geom. 68 (2004), no. 2, 185--222.

\bibitem[Z00]{z00}\textsc{K. Zuo:} On the negativity of kernels of
Kodaira-Spencer maps on Hodge bundles and applications. Asian J. of Math. 4,
(200), 270-302.
\end{thebibliography}
\end{document}